\documentclass[12pt]{article}

\makeatletter
\newfont{\footsc}{cmcsc10 at 8truept}
\newfont{\footbf}{cmbx10 at 8truept}
\newfont{\footrm}{cmr10 at 10truept}
\renewcommand{\ps@plain}{}
\makeatother 

\title{ON THE SIZE OF GRAPHS WITHOUT REPEATED CYCLE LENGTHS }

\author{ Chunhui Lai \thanks{ Project supported by  the National Science Foundation of China (No. 11401290; No.61379021),
 NSF of Fujian(2015J01018; 2014J01020; 2016J01027; 2016J01673 ), Fujian Provincial Training Foundation
for "Bai-Quan-Wan Talents Engineering", Project of Fujian Education Department (JA15295).}\\
\small School of Mathematics and Statistics, Minnan Normal University,\\
\small  Zhangzhou, Fujian 363000, P. R. of CHINA\\
\small \texttt{ laich2011@msn.cn; laichunhui@mnnu.edu.cn }\\
\small MR Subject Classifications: 05C38, 05C35\\
\small Key words: graph, cycle, number of edges}
\date{}
\begin{document}
\maketitle

\begin{abstract}
In 1975, P. Erd\"{o}s proposed the problem of determining the
maximum number $f(n)$ of edges in a  graph of $n$ vertices in
which any two cycles are of different
 lengths. In this paper, it is proved that $$f(n)\geq n+\frac{107}{3}t+\frac{7}{3}$$ for $t=1260r+169 \,\ (r\geq 1)$
 and $n \geq \frac{2119}{4}t^{2}+87978t+\frac{15957}{4}$. Consequently,
$\liminf\sb {n \to \infty} {f(n)-n \over \sqrt n} \geq \sqrt {2 +
\frac{7654}{19071}},$   which is better than the previous bounds $\sqrt
2$ [Y. Shi, Discrete Math. 71(1988), 57-71], $\sqrt {2.4}$
   [C. Lai,  Australas. J. Combin. 27(2003), 101-105].
   The conjecture $\lim_{n \rightarrow \infty} {f(n)-n\over \sqrt n}=\sqrt {2.4}$ is not true.
  \par

 \end{abstract}

\section{Introduction}
\par
Let $f(n)$ be the maximum number of edges in a graph on $n$ vertices
in which no two cycles have the same length. In 1975, Erd\"{o}s
raised the problem of determining $f(n)$ (see Bondy and Murty [1],
p.247, Problem 11).  Shi[14] proved a lower bound.
\par
\par
\noindent{\bf Theorem 1 (Shi[14])}
$$f(n)\geq n+
 [(\sqrt {8n-23} +1)/2]$$ for $n\geq 3$.
\par
 Chen, Lehel, Jacobson, and Shreve [3], Jia[5],
 Lai[6,7,8], Shi[15,17,18,19]

 obtained some results.
\par
 Boros, Caro, F\"uredi and Yuster[2]
 proved an upper bound.
\par
\noindent{\bf Theorem 2 ( Boros, Caro, F\"uredi and Yuster[2])}
  For $n$ sufficiently large, $$f(n) < n+1.98\sqrt{n}.$$
\par
Lai [9] improved the lower bound.
\par
\noindent{\bf Theorem 3 (Lai [9])}
  $$f(n)\geq n+ \sqrt {2.4} \sqrt {n}(1-o(1))$$ and proposed the following conjecture:
  \par
  \noindent{\bf Conjecture 4 (Lai [9])}
   $$\lim_{n \rightarrow \infty} {f(n)-n\over \sqrt n}=\sqrt {2.4}.$$
   \par
Lai [6] raised the following Conjecture:
 \par
\noindent{\bf Conjecture 5 (Lai [6])}
$$\liminf_{n \rightarrow \infty} {f(n)-n\over \sqrt n} \leq \sqrt {3}.$$
\par
Markstr\"om [12] raised the following problem:
 \par
\noindent{\bf Problem 6 (Markstr\"om [12])} Determine
 the maximum number of edges in a hamiltonian
graph on $n$ vertices with no repeated cycle lengths.
\par
\par
\par
 Let $f_2(n)$ be the maximum number of edges in a $2$-connected
 graph on $n$ vertices in which no two cycles have the same length. The result can be found in [2,3,14].
 \par

 The survey article on this problem can be found in Tian[20], Zhang[23], Lai and Liu [10].
\par
The progress of all 50 problems in [1] can be find in Stephen C. Locke, Unsolved problems: http://math.fau.edu/locke/Unsolved.htm
\par

A related topic concerns Entringer problem. Determine which simple graph $G$ have exactly one cycle of each
length $l,$ $3\leq l \leq v$ (see problem 10 in [1]), this problem was posed in 1973 by R. C. Entringer.
For the developments on this topic, see[4,11,12,13,16,21,22].
\par

 In this paper, we  construct a
 graph $G$ having no two cycles with the same length which leads  the following result.
  \par
\bigskip
  \noindent{\bf Theorem 7} Let $ t=1260r+169 \,\ (r\geq 1)$, then $$f(n)\geq n+\frac{107}{3}t+\frac{7}{3}$$ for
  $n\geq
  \frac{2119}{4}t^{2}+87978t+\frac{15957}{4}$. The Conjecture 4 is not true.

 \section{Proof of the theorem 7}
  {\bf Proof.} Let $t=1260r+169,r\geq 1,$
  $ n_{t}=\frac{2119}{4}t^{2}+87978t+\frac{15957}{4},$ $n\geq n_{t}.$
  We shall show that there exists a graph $G$ on $n$ vertices with $ n+\frac{107}{3}t+\frac{7}{3}$ edges such
  that all cycles in $G$ have distinct lengths.

  Now we construct the graph $G$ which consists of a number of subgraphs: $B_i$,
  ($0\leq i\leq 20t, 27t\leq i\leq 28t+64,
  29t-734\leq i\leq 29t+267, 30t-531\leq i\leq 30t+57, 31t-741\leq i\leq 31t+58,
  32t-740\leq i\leq 32t+57, 33t-741\leq i\leq 33t+57, 34t-741\leq i\leq 34t+52,
  35t-746\leq i\leq 35t+60, 36t-738\leq i\leq 36t+60,37t-738\leq i\leq 37t+799,
  i=20t+j (1\leq j \leq \frac{t-7}{6}), i=20t+\frac{t-1}{6}+j (1\leq j \leq \frac{t-7}{6}),
  i=21t+2j+1(0\leq j \leq t-1), i=21t+2j(0 \leq j \leq \frac{t-1}{2}), i=23t+2j+1
  (0 \leq j \leq \frac{t-1}{2}),$ and $i=20t+\frac{t-1}{6}$, $i=20t+\frac{t-1}{3}+\frac{t-1}{6}-1$,
  $i=20t+\frac{t-1}{3}+\frac{t-1}{6}$, $i=20t+\frac{2t-2}{3}$,
  $i=21t-2$, $i=21t-1$, $i=26t$).

  Now we define these $B_i$'s. These subgraphs all have a common vertex $x$,
  otherwise their vertex sets are pairwise disjoint.
  \par
For $1 \leq i\leq \frac{t-7}{6},$ let the subgraph $B_{20t+i}$
consist of a cycle $$xa_i^1a_i^2...a_i^{\frac{62t-8}{3}+2i}x$$  and a path:

  $$xa_{i,1}^1a_{i,1}^2...a_{i,1}^{\frac{59t-5}{6}}a_i^{\frac{61t-1}{6}+i}$$
  \par
  Based the construction,  $B_{20t+i}$ contains exactly three
cycles of lengths:
$$20t+i, 20t+\frac{t-1}{3}+i-1, 20t+\frac{2t-2}{3}+2i-1.$$
\par
   \par
For $1 \leq i\leq \frac{t-7}{6},$ let the subgraph $B_{20t+\frac{t-1}{6}+i}$
consist of a cycle $$xb_i^1b_i^2...b_i^{\frac{62t-5}{3}+2i}x$$  and a path:

  $$xb_{i,1}^1b_{i,1}^2...b_{i,1}^{10t-1}b_i^{\frac{61t-1}{6}+i}$$
  \par
  Based the construction,  $B_{20t+\frac{t-1}{6}+i}$ contains exactly three
cycles of lengths:
$$20t+\frac{t-1}{6}+i, 20t+\frac{t-1}{3}+\frac{t-1}{6}+i, 20t+\frac{2t-2}{3}+2i.$$
\par

  \par
For $0 \leq i\leq t-1,$ let the subgraph $B_{21t+2i+1}$ consist of
a cycle $$xu_i^1u_i^2...u_i^{25t+2i-1}x$$  and a path:

  $$xu_{i,1}^1u_{i,1}^2...u_{i,1}^{(19t+2i-1)/2}u_i^{(23t+2i+1)/2}$$
  \par
  Based the construction,  $B_{21t+2i+1}$ contains exactly three
cycles of lengths:
$$21t+2i+1, 23t+2i, 25t+2i.$$
\par
For $0 \leq i\leq \frac{t-3}{2},$ let the subgraph $B_{21t+2i}$
consist of a cycle $$xv_i^1v_i^2...v_i^{25t+2i}x$$  and a path:

  $$xv_{i,1}^1v_{i,1}^2...v_{i,1}^{9t+i-1}v_i^{12t+i}$$
  \par
  Based the construction,  $B_{21t+2i}$ contains exactly three
cycles of lengths:
$$21t+2i, 22t+2i+1, 25t+2i+1.$$
\par
For $0 \leq i\leq \frac{t-3}{2},$ let the subgraph $B_{23t+2i+1}$
consist of a cycle $$xw_i^1w_i^2...w_i^{26t+2i+1}x$$  and a path:

  $$xw_{i,1}^1w_{i,1}^2...w_{i,1}^{(21t+2i-1)/2}w_i^{(25t+2i+1)/2}$$
  \par
  Based the construction,  $B_{23t+2i+1}$ contains exactly three
cycles of lengths:
$$23t+2i+1, 24t+2i+2, 26t+2i+2.$$
\par

   For $58\leq i\leq t-742,$ let the subgraph  $B_{27t+i-57}$ consist of a
  cycle $$C_{27t+i-57}=xy_i^1y_i^2...y_i^{132t+11i+893}x$$  and ten paths
  sharing a common vertex $x$, the other end vertices are on the
  cycle $C_{27t+i-57}$:

  $$xy_{i,1}^1y_{i,1}^2...y_{i,1}^{(17t-1)/2}y_i^{(37t-115)/2+i}$$
  $$xy_{i,2}^1y_{i,2}^2...y_{i,2}^{(19t-1)/2}y_i^{(57t-103)/2+2i}$$
  $$xy_{i,3}^1y_{i,3}^2...y_{i,3}^{(19t-1)/2}y_i^{(77t+315)/2+3i}$$
  $$xy_{i,4}^1y_{i,4}^2...y_{i,4}^{(21t-1)/2}y_i^{(97t+313)/2+4i}$$
  $$xy_{i,5}^1y_{i,5}^2...y_{i,5}^{(21t-1)/2}y_i^{(117t+313)/2+5i}$$
  $$xy_{i,6}^1y_{i,6}^2...y_{i,6}^{(23t-1)/2}y_i^{(137t+311)/2+6i}$$
  $$xy_{i,7}^1y_{i,7}^2...y_{i,7}^{(23t-1)/2}y_i^{(157t+309)/2+7i}$$
  $$xy_{i,8}^1y_{i,8}^2...y_{i,8}^{(25t-1)/2}y_i^{(177t+297)/2+8i}$$
  $$xy_{i,9}^1y_{i,9}^2...y_{i,9}^{(25t-1)/2}y_i^{(197t+301)/2+9i}$$
  $$xy_{i,10}^1y_{i,10}^2...y_{i,10}^{(27t-1)/2}y_i^{(217t+305)/2+10i}.$$

\par
As a cycle with $d$ chords contains ${{d+2} \choose 2}$ distinct
cycles, $B_{27t+i-57}$ contains exactly 66 cycles of lengths:
 \par
 $$\begin{array}{llll}27t+i-57,& 28t+i+7,& 29t+i+210,& 30t+i,\\
  31t+i+1,& 32t+i,& 33t+i,& 34t+i-5,\\
  35t+i+3,& 36t+i+3,& 37t+i+742,& 38t+2i-51,\\
  38t+2i+216,& 40t+2i+209,& 40t+2i,& 42t+2i,\\
  42t+2i-1,& 44t+2i-6,& 44t+2i-3,& 46t+2i+5,\\
  46t+2i+744,& 48t+3i+158,& 49t+3i+215,& 50t+3i+209,\\
  51t+3i-1,& 52t+3i-1,& 53t+3i-7,& 54t+3i-4,\\
  55t+3i-1,& 56t+3i+746,& 59t+4i+157,& 59t+4i+215,\\
  61t+4i+208,& 61t+4i-2,& 63t+4i-7,& 63t+4i-5,\\
  65t+4i-2,& 65t+4i+740,& 69t+5i+157,& 70t+5i+214,\\
  71t+5i+207,& 72t+5i-8,& 73t+5i-5,& 74t+5i-3,\\
  75t+5i+739,& 80t+6i+156,& 80t+6i+213,& 82t+6i+201,\\
  82t+6i-6,& 84t+6i-3,& 84t+6i+738,& 90t+7i+155,\\
  91t+7i+207,& 92t+7i+203,& 93t+7i-4,& 94t+7i+738,\\
  101t+8i+149,& 101t+8i+209,& 103t+8i+205,& 103t+8i+737,\\
  111t+9i+151,& 112t+9i+211,& 113t+9i+946,& 122t+10i+153,\\
  122t+10i+952,& 132t+11i+894.&&\end{array}$$
 \par
 $B_{0}$ is a path with an end vertex $x$ and length $n-n_{t}$. Other $B_i$ is
 simply a cycle of length $i$.
  \par Then  $f(n)\geq n+\frac{107}{3}t+\frac{7}{3},$ for $n\geq n_{t}.$
 \par
  From the above result, we have $$\liminf_{n \rightarrow \infty} {f(n)-n\over \sqrt n}
  \geq \sqrt {2 +
\frac{7654}{19071}},$$
  which is better than the previous bounds $\sqrt 2$ (see [14]), $\sqrt {2+{2\over 5}}$
  (see [9]).
  \par
  The Conjecture 4 is not true.
  This completes
  the proof.
 \vskip 0.2in
  \par
  Combining this with Boros, Caro, F\"uredi and Yuster's upper bound, we get
  $$1.98\geq \limsup_{n \rightarrow \infty} {f(n)-n\over \sqrt n} \geq
  \liminf_{n \rightarrow \infty} {f(n)-n\over \sqrt n}\geq \sqrt {2 +
\frac{7654}{19071}}.$$
  \par

 \par

\end{document}